\documentclass{elsart3-1}

 \usepackage{graphicx}
\usepackage{amssymb}
\usepackage{amsmath}
\usepackage{setspace}

\usepackage[english,francais]{babel}

\newtheorem{e-proposition}[theorem]{Proposition}

\newtheorem{e-definition}[theorem]{Definition\rm}
\newtheorem{remark}{\it Remark\/}

\renewcommand{\theequation}{\arabic{equation}}
\newcommand{\bigo}{\mathcal{O}}
\usepackage{color}

\newcommand{\eps}{\varepsilon}

\def\og{\leavevmode\raise.3ex\hbox{$\scriptscriptstyle\langle\!\langle$~}}
\def\fg{\leavevmode\raise.3ex\hbox{~$\!\scriptscriptstyle\,\rangle\!\rangle$}}

\def\be{\begin{equation}}
\def\ee{\end{equation}}

\journal{the Acad\'emie des sciences}
\begin{document}
% place in the next line the header (rubrique) chosen for your article,
% if you know it (you can also have 2, format : Header1/Header2
\centerline{}
\begin{frontmatter}

% Title, authors and addresses

% use the thanksref command within \title, \author or \address for footnotes;
% use the ead command for the email address,
% and the form \ead[url] for the home page:
% \title{Title\thanksref{label1}}
% \thanks[label1]{}
% \author{Name\thanksref{label2}}
% \ead{email address}
% \ead[url]{home page}
% \thanks[label2]{}
% \address{Address\thanksref{label3}}
% \thanks[label3]{}
\selectlanguage{english}
\title{Asymptotic Preserving numerical schemes\\ 
for multiscale parabolic problems}

% use optional labels to link authors explicitly to addresses:
% \author[label1,label2]{}
% \address[label1]{}
% \address[label2]{}
% The [label1] can be suppressed if there is only one address for all authors

\selectlanguage{english}
\author[crouseilles,ur]{Nicolas Crouseilles},
\ead{nicolas.crouseilles@inria.fr}
%\author[ur]{H\'el\`ene Hivert}
%\ead{helene.hivert@univ-rennes1.fr}
\author[lemou,ur]{Mohammed Lemou}, 
\ead{mohammed.lemou@univ-rennes1.fr}
\author[ug]{Gilles Vilmart}
\ead{Gilles.Vilmart@unige.ch}

\address[crouseilles]{INRIA}
\address[lemou]{CNRS}
\address[ur]{IRMAR, Universit\'e de Rennes 1. Campus de Beaulieu. 35000 Rennes, France.}
\address[ug]{Universit\'e de Gen\`eve, Section de math\'ematiques, 2-4 rue du Li\`evre, CP 64, 1211 Gen\`eve 4, Switzerland.}

% If you know the dates of reception, and acceptation you can put them now;
%  idem the name of the person presenting the Note

\medskip
\begin{center}
{\small Received *****; accepted after revision +++++\\
Presented by £££££}
\end{center}

\begin{abstract}\selectlanguage{english}%
We consider a class of multiscale parabolic problems with diffusion coefficients 
oscillating in space at a possibly small scale $\varepsilon$.
Numerical homogenization methods are popular for such problems, because they
capture efficiently the asymptotic behaviour as $\varepsilon \rightarrow 0$, without using a dramatically fine spatial discretization at the scale of the fast oscillations.
However, known such homogenization schemes are in general not accurate for both the highly oscillatory regime $\varepsilon \rightarrow 0$ and the non oscillatory regime $\varepsilon \sim 1$.
In this paper, we introduce an Asymptotic Preserving method based on an exact micro-macro decomposition of the solution which remains consistent for both regimes.

%{\it To cite this article: N.
%Crouseilles, M. Lemou, G. Vilmart, C. R. Acad. Sci. Paris, Ser. I 340 (2015) ????.}
\vskip 0.5\baselineskip

\selectlanguage{francais}
% Text of abstract in French

\noindent{\bf R\'esum\'e} \vskip 0.5\baselineskip 
\noindent
{\bf Sch\'emas num\'eriques Asymptotic Preserving pour les probl\`emes paraboliques  multi-\'echelles. }
On consid\`ere une classe de probl\`emes paraboliques multi-\'echelles dont les coefficients de diffusion
oscillent rapidement en espace \`a une \'echelle $\varepsilon$ possiblement petite.
Les m\'ethodes num\'eriques d'homog\'en\'eisation sont populaires pour ces probl\`emes, car elles
capturent efficacement le comportement asymptotique lorsque $\varepsilon \rightarrow 0$,
sans utiliser une discr\'etisation spatiale aussi fine que l'\'echelle des oscillations rapides, 
comme le n\'ecessiteraient les m\'ethodes non-raides standards.
Cependant, les sch\'emas d'homog\'en\'eisation existants ne sont en g\'en\'eral pas pr\'ecis dans les deux r\'egimes oscillant  $\varepsilon \rightarrow 0$ et non-oscillant $\varepsilon \sim 1$.
Dans ce travail, nous introduisons une m\'ethode Asymptotic Preserving bas\'ee sur une d\'ecomposition micro-macro exacte qui reste consistante pour les deux r\'egimes.

%{\it Pour citer cet article~: N. Crouseilles, M. Lemou,  G. Vilmart, C. R. Acad. Sci.
%Paris, Ser. I 340 (2005). ??????}
\end{abstract}
\end{frontmatter}

\selectlanguage{francais}
\section*{Version fran\c{c}aise abr\'eg\'ee}

L'objectif est de construire des sch\'emas num\'eriques pour des probl\`emes paraboliques \eqref{eq:1}
o\`u les coefficients de diffusion sont hautement oscillants en espace. Ce type de mod\`ele intervient dans des probl\`emes de propagation en milieux exhibants une structure p\'eriodique. Dans de nombreux cas, la taille de la p\'eriode est petite par rapport 
\`a la taille caract\'eristique du milieu, et en notant par $\varepsilon$ leur rapport,  une analyse asymptotique est n\'ecessaire 
pour \'etudier le comportement de la solution quand $\varepsilon\rightarrow 0$. 
Cette analyse a \'et\'e conduite dans \cite{BLP78,BFM92,All92,JKO94} \`a l'aide de la convergence double-\'echelle. 

D'un point de vue num\'erique, une approche directe s'av\`ere tr\`es co\^uteuse puisque les param\`etres 
num\'eriques des maillages doivent r\'esoudre la plus petite \'echelle $\varepsilon$.
Les m\'ethodes num\'eriques d'homog\'en\'eisation comme la ``heterogeneous multiscale method'' (HMM) \cite {EE03} (voir \cite{AbV12} dans le context de probl\`emes paraboliques lin\'eaires, et l'article de revue \cite{AEE12}) permettent de calculer efficacement la solution homog\'en\'eis\'ee $u^0$ mais aussi la solution oscillante $u^\eps$ \`a $\eps$ fix\'e 
en se basant sur l'approximation du mod\`ele asymptotique, qui suppose que $\varepsilon$ est tr\`es petit, ainsi que des techniques de correcteurs \cite{BFM92,JKO94}. Cependant, si 
$\varepsilon$ n'est pas petit, ce type d'approche tombe en d\'efaut. 
%{\color{red} On mentionne \'egalement la ``multiscale finite element method'' (msFEM) \cite{HoX97} (voir l'article de revue \cite{EfH09}), qui permet de calculer directement la solution  $u^\varepsilon$ \`a l'aide d'une base adapt\'ee de fonctions oscillantes, mais dont le co\^ut augment quand $\epsilon \rightarrow 0$.}
   
Dans ce travail, nous proposons une m\'ethode qui permet d'approcher la solution $u^\eps$ pour n'importe quelle valeur de $\varepsilon\in (0, 1]$, % fix\'ee, 
\`a param\`etres num\'eriques fix\'es ind\'ependament de $\eps$.
Ce type d'approche est dit "Asymptotic Preserving" \cite{Jin99} : un tel sch\'ema 
est consistant avec le probl\`eme initial pour tout $\varepsilon$ fix\'e et d\'eg\'en\`ere 
quand $\varepsilon\rightarrow 0$ en un sch\'ema num\'erique consistant avec le mod\`ele asymptotique. 

Notre approche est bas\'ee sur une reformulation du probl\`eme initial en un probl\`eme augment\'e (voir \eqref{eq:2} 
et les notations \eqref{notations}) satisfait par $U^\eps(t,x,y)$, 
dans lequel les \'echelles lentes $x$ et  rapides $y=x/\varepsilon$ sont consid\'er\'ees ind\'ependantes. 
La solution $u^\eps(t,x)$ 
du probl\`eme initial peut alors \^etre retrouv\'ee gr\^ace \`a la relation $U^\eps(t,x, y=x/\varepsilon)=u^\eps(t,x)$. 
Dans le probl\`eme \eqref{eq:2}, une raideur appara\^it devant le terme $LU^\eps=\nabla_y\cdot (a(x, y)\nabla_y U^\eps)$ 
qui permet de s'inspirer des m\'ethodes de d\'eveloppement asymptotique largement 
utilis\'ee  en th\'eorie cin\'etique pour construire 
des sch\'emas num\'eriques multi-\'echelles   (voir \cite{LeM08}). 
Nous d\'ecomposons en effet la solution du probl\`eme augment\'e $U^\eps$ sous la forme $U^\eps(t,x,y)=F^\eps(t,x)+G^\eps(t,x,y)$, 
o\`u $F^\eps$ est la projection orthogonale de $U^\eps$ sur le noyau de $L$. 
%Cette approche %d\'ecomposition 
%est analogue aux d\'ecompositions micro-macro utilis\'ees en th\'eorie cin\'etique. 
Cette d\'ecomposition permet de reformuler de fa\c con \'equivalente 
le probl\`eme augment\'e en un syst\`eme d'\'equations micro-macro satisfait par $F^\eps$ et $G^\eps$. 
Notons que ce syst\`eme ne contient aucune approximation et reste exact pour toute 
valeur de $\varepsilon$. 

Nous nous basons ensuite sur cette d\'ecomposition pour construire un sch\'ema num\'erique 
multi-\'echelles. Gr\^ace \`a une 
discr\'etisation  en temps semi-implicite (inspir\'ee de \cite{LeM08}), un sch\'ema Asymptotic Preserving 
est alors obtenu. Ce sch\'ema n\'ecessite l'inversion au niveau num\'erique de l'op\'erateur $L$ (\`a $x$ fix\'e), 
comme dans le cas de la r\'esolution num\'erique  du probl\`eme homog\'en\'eis\'e (voir \cite{EE03}). 

Bien que la m\'ethodologie soit pr\'esent\'ee ici en dimension quelconque, nous effectuons 
des tests num\'eriques uniquement en dimension $1$.  Le bon comportement du sch\'ema 
pour diff\'erents $\varepsilon\in (0, 1]$ est mis en \'evidence en le comparant \`a une 
m\'ethode directe et \`a la solution du probl\`eme homog\'en\'eis\'e. L'analyse
du cas multi-dimensionnel avec une m\'ethode d'\'el\'ements finis est en cours d'\'etude. Cette note est une version abr\'eg\'ee du travail plus d\'etaill\'e   \cite{CLV2}.

% Text of your Version française abr\'eg\'ee here.
% Note you do not need to repeat here equations that you use in the
% main text - for example 'voir (3)' is quite acceptable.

\selectlanguage{english}

\vspace{-5mm}

\section{Introduction}
For $T>0$ and a smooth bounded domain $\Omega \subset \mathbb{R}^d$, we consider the following class of parabolic problems
\be
\label{eq:1}
\partial_t u^\varepsilon= \nabla_x \cdot  \left[a(x,x/\varepsilon) \nabla_x u^\varepsilon \right] + f(t,x), \;\; t\in (0, T), \;\; x\in \Omega, 
\ee
where $u^\eps(t=0, x)=g(x)$ is a given initial condition in $L^2(\Omega)$, $f \in L^2(0,T,L^2(\Omega))$ is a given source term, and we take for simplicity homogeneous Dirichlet boundary conditions $u^\varepsilon=0$ in $(0,T)\times \partial\Omega$.
The tensor $a(x,y) \in \mathbb{R}^{d\times d}$ 
is assumed symmetric, uniformly elliptic and bounded, and periodic with respect to the variable $y=x/\varepsilon \in Y=(0,1)^d$. 
The homogenization analysis as $\varepsilon\rightarrow 0$ of such a multiscale problem is well-known and can be done using a two-scale convergence analysis, see \cite{BLP78,BFM92,All92,JKO94}.
Problem \eqref{eq:1} admits a unique solution $u^\varepsilon$ is the space $L^2(0,T; H_0^1(\Omega))$ which converges towards an asymptotic solution $u^0$ as $\eps\rightarrow 0$, 
$$
u^{\varepsilon}\rightarrow u^0 \hbox{ strongly in } L^{2}(0,T;L^2(\Omega)),\quad
u^{\varepsilon}\rightharpoonup u^0 \hbox{ weakly in } L^{2}(0,T;H_0^1(\Omega)),
$$
where $u^0 \in L^2(0,T; H_0^1(\Omega))$ solves an effective non-oscillatory problem of the same form as \eqref{eq:1}, 
\begin{eqnarray}
\label{eq-allaire-1}
-\nabla_y \cdot[a(x,y) [\nabla_x u^0 + \nabla_y u_1]] &=& 0, \quad x\in \Omega,y\in Y,\\
\label{eq-allaire-2}
\partial_t u^0-\nabla_x \cdot\left[  \int_Y a(x,y) [\nabla_x u^0 + \nabla_y u_1] dy \right] &=&  f, \quad x\in \Omega,  
\end{eqnarray}
 with the same initial and boundary conditions for $u^0$ in \eqref{eq-allaire-2} as for $u^\eps$ and 
involving the elliptic ``cell problem'' \eqref{eq-allaire-1} 
with solution $u_1(t,x,\cdot) \in H^1_{per}(Y)$, periodic with zero average with respect to the second variable $y$.
Taking advantage of the separation of micro and macro scales, numerical homogenization methods exploit the above homogenization result to compute efficiently the solution of \eqref{eq:1} in the asymptotic regime $\varepsilon \rightarrow 0$.
For instance, such an efficient method is the heterogeneous Multiscale Method (HMM) \cite{EE03} (see also the review \cite{AEE12}) which relies on a coupling of micro and macro finite element methods applied to \eqref{eq-allaire-1}-\eqref{eq-allaire-2} 
(see \cite{AbV12} in the context of parabolic multiscale problems \eqref{eq:1}). 
Having a computational cost independent of the smallness of $\varepsilon$, HMM permits to approximate not only the asymptotic solution $u^0$, but also the oscillatory solution $u^\varepsilon$ and its gradient $\nabla_x u^\eps$ for a fixed small $\varepsilon$ using the approximation 
$u^\varepsilon(t,x) \simeq u^0(t,x) + \varepsilon u_1(t,x, y=x/\varepsilon)$  
based on corrector techniques. 
This latter approximation of $u^\varepsilon$ is consistent only for small values of $\varepsilon$ but not for $\eps$ close to $1$.
We also mention the ``multiscale finite element method'' (msFEM) \cite{HoX97} (see the review \cite{EfH09}), which permits to compute the oscillatory solution $u^\varepsilon$ by using an enriched finite element space, but where the computational cost grows as $\epsilon \rightarrow 0$.
The aim of this paper is to introduce a micro-macro decomposition which permits to approximate $u^\varepsilon$ accurately for both regimes $\varepsilon \rightarrow 0$  or $\varepsilon\sim 1$ at a cost independent of $\varepsilon$, in the spirit of Asymptotic Preserving schemes \cite{Jin99}.\looseness-3

\medskip 
\begin{remark} \label{rem:1}
The asymptotic parabolic problem \eqref{eq-allaire-2} is non-stiff and can be written in the form $\partial_t u^0= \overline{D}u^0 + f$ 
where $\overline{D} \phi = \nabla_x \cdot (a^0 \nabla_x \phi)$ and $a^0$ is the so-called homogenized tensor. In other words, the asymptotic problem has the same form as \eqref{eq:1} with $a^\varepsilon$ replaced by $a^0$.
In dimension $d=1$, the homogenized tensor is given by the harmonic average 
$a^0(x)= ({\int_Y a(x,y)^{-1}dy})^{-1}$. However, in multiple dimensions, there is no such a simple formula.
The calculation of the asymptotic diffusion coefficient $a^0(x)$ at a point $x$ of space then requires the resolution of 
an elliptic type problem like \eqref{eq-allaire-1}. We refer to the review \cite{AEE12}. % in the context HMM for details.
\end{remark}

\vspace{-0.5cm}

\section{Exact micro-macro decomposition}
In this section, we introduce a micro-macro decomposition of the oscillatory solution $u^\varepsilon$ of \eqref{eq:1} which remains exact for all $\varepsilon \in (0,1)$ and we study its behavior for $\varepsilon \rightarrow 0$. 
We emphasize that our analysis is only formal, and a rigorous study of the approach is currently under investigation \cite{CLV2}. 
In the spirit of two-scale convergence analysis, we introduce a function $U^\varepsilon: (0, T)\times \Omega\times Y \rightarrow \mathbb{R}$,  periodic with respect to the third variable $y\in Y=(0, 1)^d$, and 
such that $U^\varepsilon$ coincides with $u^\varepsilon$, the solution to \eqref{eq:1}, on the diagonal, {\it i.e.} 
\be
\label{diagonal}
U^\varepsilon(t,x,x/\varepsilon)= u^\varepsilon(t,x),
\ee
and we obtain that $U^\varepsilon$ solves the following augmented problem  
\be
\label{eq:2}
\partial_t U^\varepsilon =\frac{1}{\varepsilon^2} LU^\varepsilon+\frac{1}{\varepsilon}BU^\varepsilon + DU^\varepsilon + f. 
\ee
Here, we use the following notations for all functions $\phi$,
\be
\label{notations}
L\phi =  \nabla_y \cdot [a(x,y)\nabla_y \phi], \;\;\; D \phi = \nabla_x\cdot[a(x,y) \nabla_x \phi], \;\; B \phi =  \nabla_x \cdot[a(x,y)\nabla_y \phi] + \nabla_y\cdot [a(x,y)\nabla_x \phi]. 
\ee
We then choose appropriate initial and boundary conditions for \eqref{eq:2},
\be \label{eq:2cond}
U^\eps(0, x, y)=g(x) \mbox{ for } x\in \Omega, y\in Y, \quad
U^\eps(t, x, y)=\eps\left(u_1(t,x,y)-u_1(t,x,\frac x \eps)\right) \mbox{ for } t\in (0, T), x \in \partial\Omega, y\in Y,
\ee
where $u_1$ is given by \eqref{eq-allaire-1}.
Consider now the linear operator $L$ in \eqref{notations} defined for a fixed $x$ on $H^1_{per}(Y)$, the space of periodic functions in $H^1(Y)$; this operator is self-adjoint with respect to the $L^2(Y)$ scalar product. Its kernel 
is the set of constant functions (with respect to $y$) and the $L^2$ orthogonal projector on this kernel is the average projection operator $\Pi \phi := \int_Y \phi(y) dy$.
Moreover, $L$ is an isomorphism between the Hilbert space 
$W_{per}(Y)=\{\phi \in H^1_{per}(Y)\ ;\ \Pi \phi=0\}$ and its dual $(W_{per}(Y))'$. 
Following \cite{LeM08} in the context of kinetic theory, 
we now perform an exact micro-macro decomposition of $U^\varepsilon$ by setting $F^\varepsilon = \Pi U^\varepsilon$,
$G^\varepsilon = (I-\Pi) U^\varepsilon$, where $I$ is the identity operator
\be
\label{U}
U^\varepsilon(t, x,y) = F^\varepsilon(t, x) + G^\varepsilon(t, x,y). 
\ee
Inserting this decomposition into \eqref{eq:2}, 
and applying respectively $\Pi$ and $(I-\Pi)$ 
leads to the following model for the micro-macro decomposition of $u^\eps(t,x)=F^\eps(t,x)+G^\eps(t,x,x/\eps)$,
\begin{eqnarray}
\label{eq:5}
\partial_t F^\varepsilon &=& \frac{1}{\varepsilon}\Pi BG^\varepsilon + \Pi DF^\varepsilon+\Pi DG^\varepsilon + f, \\
&&F^\eps(0,x)=g(x),x\in \Omega,\quad F^\eps(t,x)=-\eps u_1(t,x,\frac x\eps),\quad t\in(0,T),x\in\partial \Omega, \nonumber\\  
\label{eq:6}
\partial_t G^\varepsilon &=&\frac{1}{\varepsilon^2}LG^\varepsilon
+ \frac{1}{\varepsilon}(I-\Pi)B[F^\varepsilon+G^\varepsilon] + (I-\Pi)D(F^\varepsilon+G^\varepsilon), \\
&&G^\eps(0,x,y)=0,x\in \Omega,y\in Y,\quad G^\eps(t,x,y)= \eps u_1(t,x, y),\quad t\in(0,T),x\in\partial \Omega,y\in Y, \nonumber
\end{eqnarray}
where $u_1$ is given in \eqref{eq-allaire-1}, and we note that $\Pi B\phi=0$ for all $\phi$ independent of $y$.

We now study formally the asymptotic limits of $G^\eps$ and $F^\eps\rightarrow \overline F$ as $\varepsilon \rightarrow 0$.
Multiplying both sides of \eqref{eq:6} by $\varepsilon$, using $\partial_t G^\varepsilon=\bigo(1)$ and setting formally $\varepsilon\rightarrow 0$, we deduce
\be
\label{eq:62}
\varepsilon^{-1} G^\varepsilon \rightarrow -L^{-1} B \overline F. 
\ee
Injecting \eqref{eq:62} into \eqref{eq:5} we deduce formally the asymptotic limit as $\varepsilon \rightarrow 0$,
\be
\label{eq:7}
\partial_t \overline  F = \overline{D} \,\overline F +f, \mbox{ with }  \overline{D}\phi := -\Pi B L^{-1} B\phi + \Pi D \phi. 
\ee
Comparing with the homogenized problem \eqref{eq-allaire-1}-\eqref{eq-allaire-2},
we emphasize that the asymptotic limits as $\varepsilon \rightarrow 0$  in \eqref{eq:62} and \eqref{eq:7} coincide respectively with $u^0$ and $u_1$ in \eqref{eq-allaire-1},\eqref{eq-allaire-2}.
Precisely, we have $F^\varepsilon \rightarrow \overline F = u^0$ and $\varepsilon^{-1} G^\varepsilon \rightarrow -L^{-1}Bu^0 = u_1$ for $\eps\rightarrow 0$, and we see that the parabolic effective problem \eqref{eq-allaire-2} is equivalent to \eqref{eq:7}. 
Therefore, the decomposition \eqref{diagonal}-\eqref{U} of $u^\varepsilon(t, x)$
into the non-oscillatory part $F^\varepsilon(t,x)$ and the oscillatory part $G^\varepsilon(t,x,x/\varepsilon)$ can indeed be interpreted as a generalization for $\varepsilon\in(0,1)$ of the approximation $u^\varepsilon \simeq u^0 + \varepsilon u_1$ valid only for small values of  $\varepsilon$.
\begin{remark} \label{rem:2}
Notice that using the simpler homogeneous boundary condition $U(t,x,y)=0$ for $t\in (0, T), x \in \partial\Omega, y\in Y$ 
(instead of \eqref{eq:2cond}) would yield an undesired boundary layer. Indeed, it is a classical issue for corrector techniques (see \cite{JKO94}) that the limit $u_1(x,y)$ appearing in \eqref{eq:62} with $y=x/\eps$ does not satisfy homogeneous Dirichlet boundary
conditions for $x\in \partial \Omega$.
\end{remark} 

\vspace{-0.5cm}

\section{Asymptotic Preserving numerical method}

The main goal of this section is to propose an Asymptotic Preserving numerical method for \eqref{eq:5}-\eqref{eq:6}, used to approximate the solution $u^\varepsilon$ of \eqref{eq:1}, and with a computational cost independent of $\eps\in(0,1)$. 
We focus on the time discretization, considering a mesh of the time interval $[0, T]$: $t^n=n\Delta t$, 
with $n\in \mathbb{N}$ and $\Delta t$ the time step. Then, we denote by $F^n$ (resp. $G^n$) 
an approximation of $F^\eps(t^n)$ (resp. $G^\eps(t^n)$). 

\noindent The stiffest term in \eqref{eq:6} has to be considered implicit to ensure stability as $\varepsilon\rightarrow 0$ 
whereas all the other terms can be considered explicit. Then, a natural first order time discretization of \eqref{eq:6}
is (see \cite{LeM08})
\be
\label{micro-scheme1}
 G^{n+1} =\left(I-\frac{\Delta t}{\varepsilon^2}L\right)^{-1}\left[G^n + \frac{\Delta t}{\varepsilon}(I-\Pi)(B+\varepsilon D)(F^n+G^n)\right]. 
\ee
Now, to recover the correct asymptotic behavior, the time discretization of \eqref{eq:5} is 
\be
\label{macro-scheme1}
F^{n+1} =F^n+ \frac{\Delta t}{\varepsilon}\Pi B G^{n+1}  + \Delta t\Pi D(F^n+G^{n+1}) + f.  
\ee
To make explicitly appear the asymptotic model, we now propose a suitable transformation 
of the macro part \eqref{macro-scheme1} (following \cite{Lem10}). To do that, we consider the 
following Duhamel formula for \eqref{eq:6}
$$
\partial_t (e^{-tL/\varepsilon^2} G^\eps) = e^{-tL/\varepsilon^2}  \left[ \frac{1}{\varepsilon}(I-\Pi)(BF^\eps+BG^\eps)+(I-\Pi)(DF^\eps+DG^\eps) \right]. 
$$
Integrating between $t^n$ and $t^{n+1}$ and performing some first order (in time) approximation leads to 
\begin{eqnarray}
G^{n+1} &\approx& e^{\Delta t L/\varepsilon^2} G^n +  \frac{1}{\varepsilon} \int_{t^n}^{t^{n+1}} e^{(t^{n+1}-s)L/\varepsilon^2}ds (I-\Pi)BF^n\nonumber\\
\label{eqGtmp}
&&+  \frac{1}{\varepsilon} \int_{t^n}^{t^{n+1}} e^{(t^{n+1}-s)L/\varepsilon^2}ds (I-\Pi)BG^n + \int_{t^n}^{t^{n+1}} e^{(t^{n+1}-s)L/\varepsilon^2}ds (I-\Pi)(DF^n+DG^n). 
\end{eqnarray}
Now, our goal is to derive an approximation of $G^{n+1}$ 
which will be inserted in the right hand side of \eqref{macro-scheme1}. This means that,  
 in the approximation of $G^{n+1}$, any term of the order of $\Delta t$ can be neglected if needed. 

In this spirit, the time integrals in \eqref{eqGtmp} are approximated as follows. 
The first one is calculated exactly to get 
$\int_{t^n}^{t^{n+1}} e^{(t^{n+1}-s)L/\varepsilon^2}ds = -\varepsilon^2 (1-e^{\Delta t L/\varepsilon^2})L^{-1}$ 
whereas the second and the third ones are approximated using a midpoint formula 
$\int_{t^n}^{t^{n+1}} e^{(t^{n+1}-s)L/\varepsilon^2}ds\approx \Delta t e^{\Delta t L/(2\varepsilon^2)}$. 
Additionally, for non small $\varepsilon$, up to terms of order $\Delta t$, $e^{\Delta tL/\varepsilon^2}$ can be replaced 
by $e^{-\Delta t/\varepsilon^2}$; for small $\varepsilon$, both $e^{\Delta tL/\varepsilon^2}$ and $e^{-\Delta t/\varepsilon^2}$ 
go to zero. Therefore, $e^{\Delta tL/\varepsilon^2}$ may be replaced by $e^{-\Delta t/\varepsilon^2}$. 
Let us emphasize that all the approximations are consistent in time with the continuous model. 
Finally, we get 
$$
G^{n+1} = e^{-\Delta t/\varepsilon^2} G^n  -\varepsilon (1-e^{-\Delta t /\varepsilon^2})L^{-1}  (I-\Pi)BF^n 
+ \frac{\Delta t}{\varepsilon} e^{-\Delta t/(2\varepsilon^2)}(I-\Pi)BG^n+\Delta t e^{-\Delta t/(2\varepsilon^2)} (I-\Pi)(DF^n+DG^n). 
$$
%\begin{eqnarray*}
%G^{n+1} &=& e^{-\Delta t/\varepsilon^2} G^n  -\varepsilon (1-e^{-\Delta t /\varepsilon^2})L^{-1}  (I-\Pi)BF^n \nonumber\\
%&&+ \frac{\Delta t}{\varepsilon} e^{-\Delta t/(2\varepsilon^2)}(I-\Pi)BG^n+\Delta t e^{-\Delta t/(2\varepsilon^2)} (I-\Pi)(DF^n+DG^n). 
%\end{eqnarray*}
This expression of $G^{n+1}$ is injected in \eqref{macro-scheme1} and using $-\Pi BL^{-1}(I-\Pi)BF^n + \Pi DF^n = \overline{D}F^n$, we have 
\be
\label{macro-man}
F^{n+1} = F^n +\Delta t(1-e^{-\Delta t/\varepsilon^2})\overline{D}F^n+\Delta t e^{-\Delta t/\varepsilon^2}\Pi DF^n+\frac{ \Delta t e^{-\Delta t/\varepsilon^2} }{\varepsilon}\Pi B  G^n + \Delta t \Pi DG^{n+1}+ \Delta t f. 
%&+&  \frac{\Delta t^2}{\varepsilon^2} e^{-\Delta t /(2\varepsilon^2)}\Pi B  (I-\Pi)BG^n \nonumber\\
%&+& \frac{ \Delta t^2}{\varepsilon} e^{-\Delta t /(2\varepsilon^2)} \Pi B(I-\Pi)(DF^n+DG^n)  + \Delta t \Pi DG^{n+1} + \Delta t S. 
\ee
The numerical scheme \eqref{micro-scheme1}-\eqref{macro-man} enjoys the Asymptotic Preserving property: 
{\it (i)} for a fixed $\varepsilon>0$, it is a first order approximation of \eqref{eq:5}-\eqref{eq:6};  
{\it (ii)} for a fixed $\Delta t$,  \eqref{micro-scheme1}-\eqref{macro-man} 
is uniformly stable with respect to $\varepsilon$ and degenerates into a consistent discretization 
of the asymptotic model \eqref{eq:7}. 

\noindent Note that up to first order terms in $\Delta t$, an implicit time discretization for the term ${\overline{D}} F$ could also 
be considered (or alternatively an explicit stabilized method as in \cite{AbV12}), which enables to avoid the parabolic CFL condition in the asymptotic regime.  Let us also remark that the numerical scheme  \eqref{micro-scheme1}-\eqref{macro-man} only requires the inversion 
of elliptic type operator in dimension $d$, which is also needed for the numerical resolution of the asymptotic model 
\eqref{eq:7}. Hence, the additional computation induced by our approach only lies in the numerical 
approximation of $\Pi, B, D$.

\section{Numerical results}

We consider a simpler case to illustrate the efficiency of our approach, considering \eqref{eq:1} in dimension $d=1$
together with the following diffusion coefficient $A(x/\varepsilon) = 1.1+\sin(x/\varepsilon)$ 
and a null right hand side $f(x)=0$ 
with the initial condition $g(x)=\sin(2\pi x)$, $x\in [0, 1]$. 
We then compare three different approaches for approximating the solution $u^\varepsilon$ to \eqref{eq:1}: 
{\it (i)} a direct approach ('REF') based on a first order explicit time integrator of \eqref{eq:1} 
in which the mesh parameters are adapted to the smallness of $\varepsilon$ (this will serve as a reference 
for comparison, obviously, another time integrator can be used, such as an implicit one);   
{\it (ii)} the exact micro-macro decomposition approach ('EMM') based on \eqref{micro-scheme1}-\eqref{macro-man} 
and {\it (iii)} the asymptotic model ('HMM') $u^\eps\simeq u^0 + \eps u_1$ with $u^0,u_1$ given by \eqref{eq-allaire-1}-\eqref{eq-allaire-2}. 
The spatial discretization (in $x$ and $y$ directions) is performed 
using a standard second order finite volume method. 
We denote by $N_x$ (resp. $N_y$) the number of points in the $x$ (resp. $y$) direction, and 
$\Delta x=1/N_x$, $\Delta y=1/N_y$ are the mesh size in $x$ and $y$ directions. 
For REF, we choose $\Delta t=0.05\Delta x^2$  
where $\Delta x$ will be chosen small enough to capture the oscillations of size $\varepsilon$. 
For EMM and HMM, we choose $\Delta t = 0.2\Delta x^2$.  
Once $U^\eps(t, x, y)$ is computed on the discrete meshes, we need to interpolate it at $y=x/\varepsilon$, $y$ being a periodic variable; 
this is done using trigonometric interpolation which ensures spectral accuracy.  
As a diagnostic for EMM, we consider a reconstructed solution on a refined mesh by using the following 
linear interpolation for $x\in [x_i, x_{i+1}],  i=0, \dots, N_x-1$,
\be
\label{ureconstruct}
{u}^\eps(t^n,x)\approx  \frac{(x_{i+1}-x)}{\Delta x} \left(F^n(x_i)+G^n(x_i, \frac x\varepsilon)\right) + \frac{(x-x_{i})}{\Delta x} \left(F^n(x_{i+1})+G^n(x_{i+1}, \frac x\varepsilon)\right),    
\ee
where $(F^n, G^n)$ given by \eqref{macro-man}-\eqref{micro-scheme1}, which enables to recover the small-scale information. For HMM, we use the same reconstruction \eqref{ureconstruct} 
with $u^0$ and $\eps u_1$ at $t=t_n$ (given in \eqref{eq-allaire-1}-\eqref{eq-allaire-2}) instead of $F^n$ and $G^n$, respectively. 
Hence, this enables 
to have an approximation of the EMM and HMM solutions as well as their spatial derivative 
on a refined mesh on which the reference solution has been obtained.

\begin{figure}[t]
\vspace{-1.cm}
\begin{tabular}{cccc}
\vspace{-2.cm}
&$\varepsilon=1$&
$\varepsilon=0.1$&
$\varepsilon=0.01$\\
\begin{minipage}{1ex}\vspace{-7.5cm} \rotatebox[origin=lB]{90}{solutions for $u^\eps$}\end{minipage}&
\vspace{-2.cm}\includegraphics[width=5.3cm]{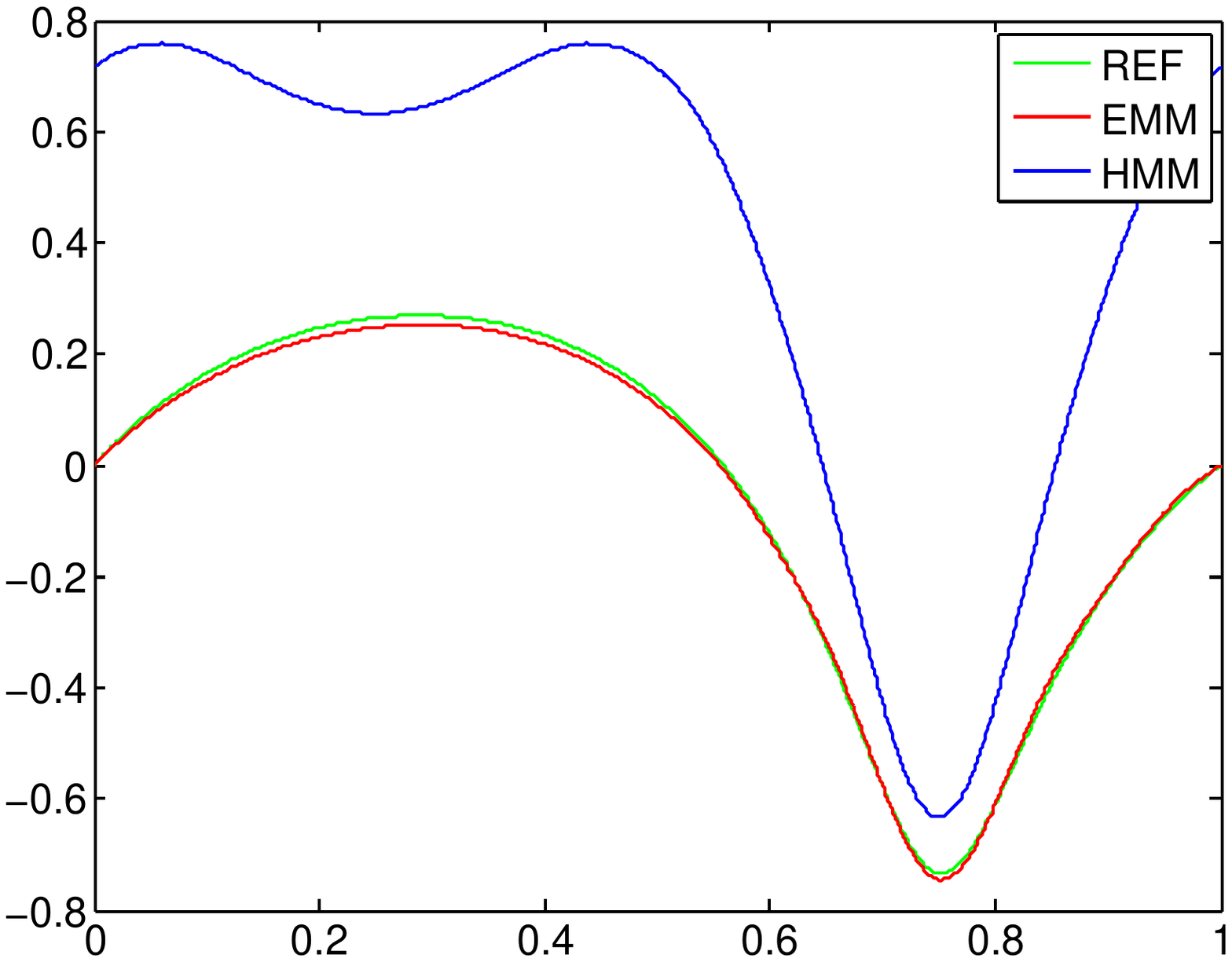} &
\includegraphics[width=5.3cm]{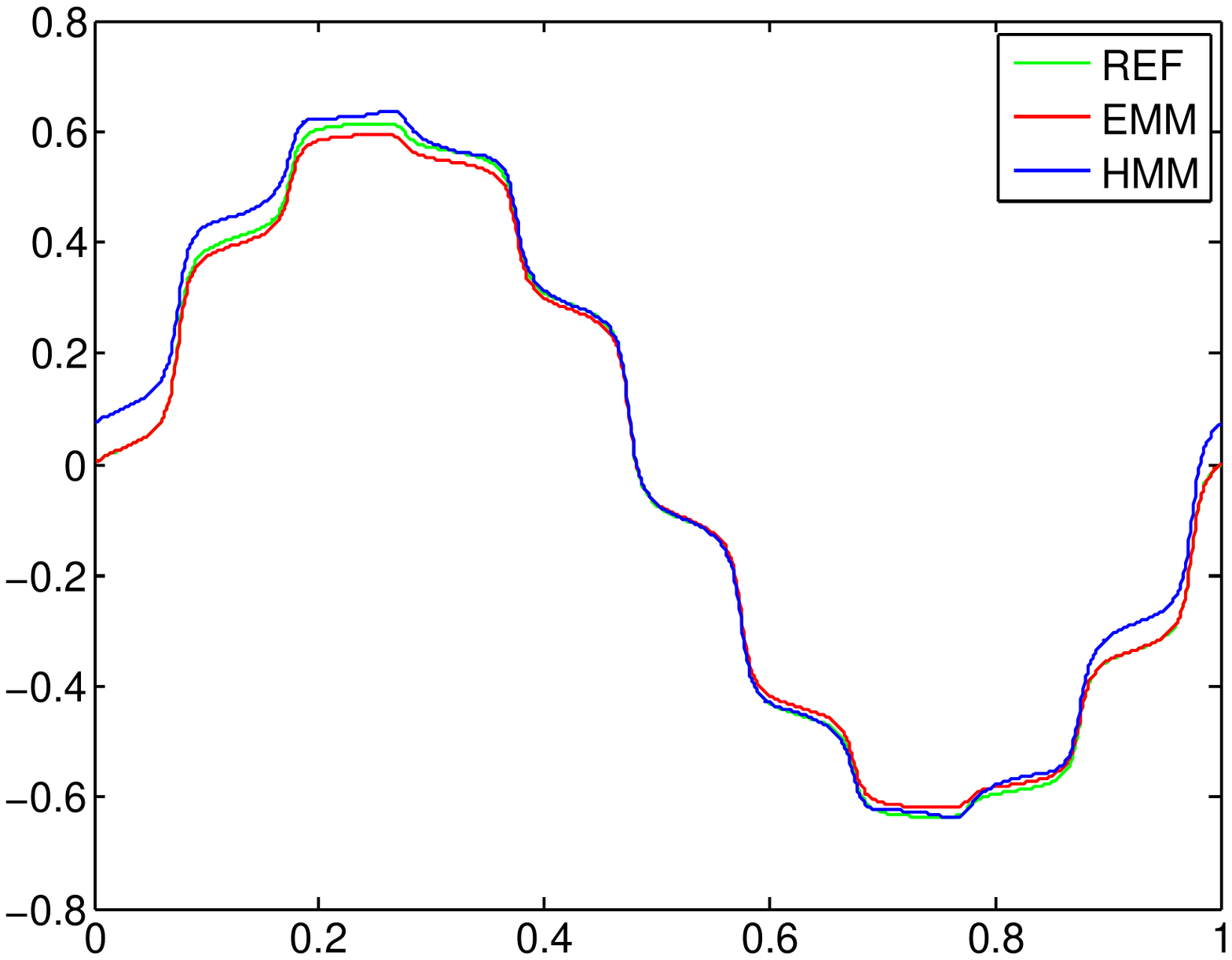} &
\vspace{-2.cm}\includegraphics[width=5.3cm]{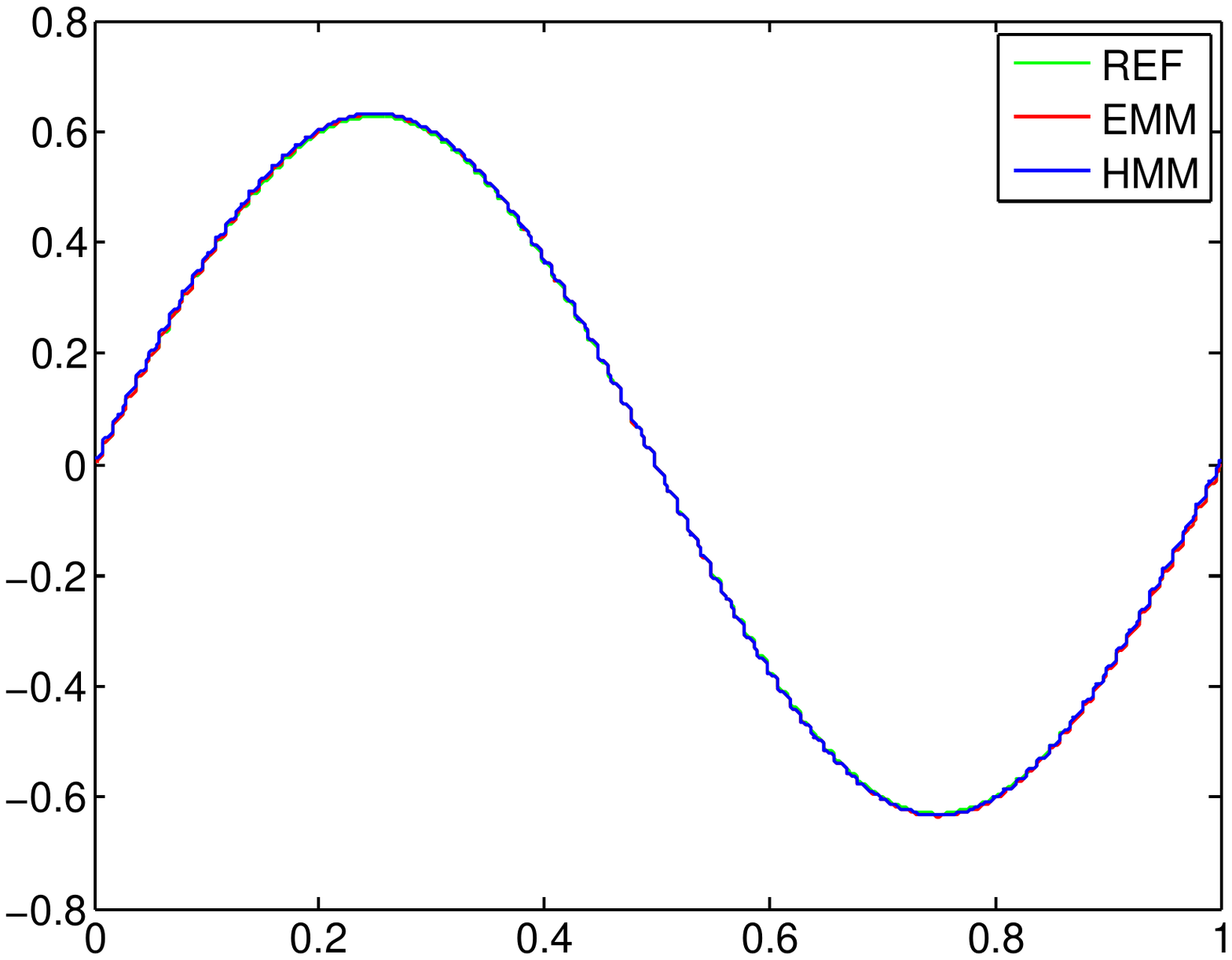} \\
\begin{minipage}{1ex}\vspace{-7.5cm} \rotatebox[origin=lB]{90}{errors for $u^\eps$}\end{minipage}&
\includegraphics[width=5.3cm]{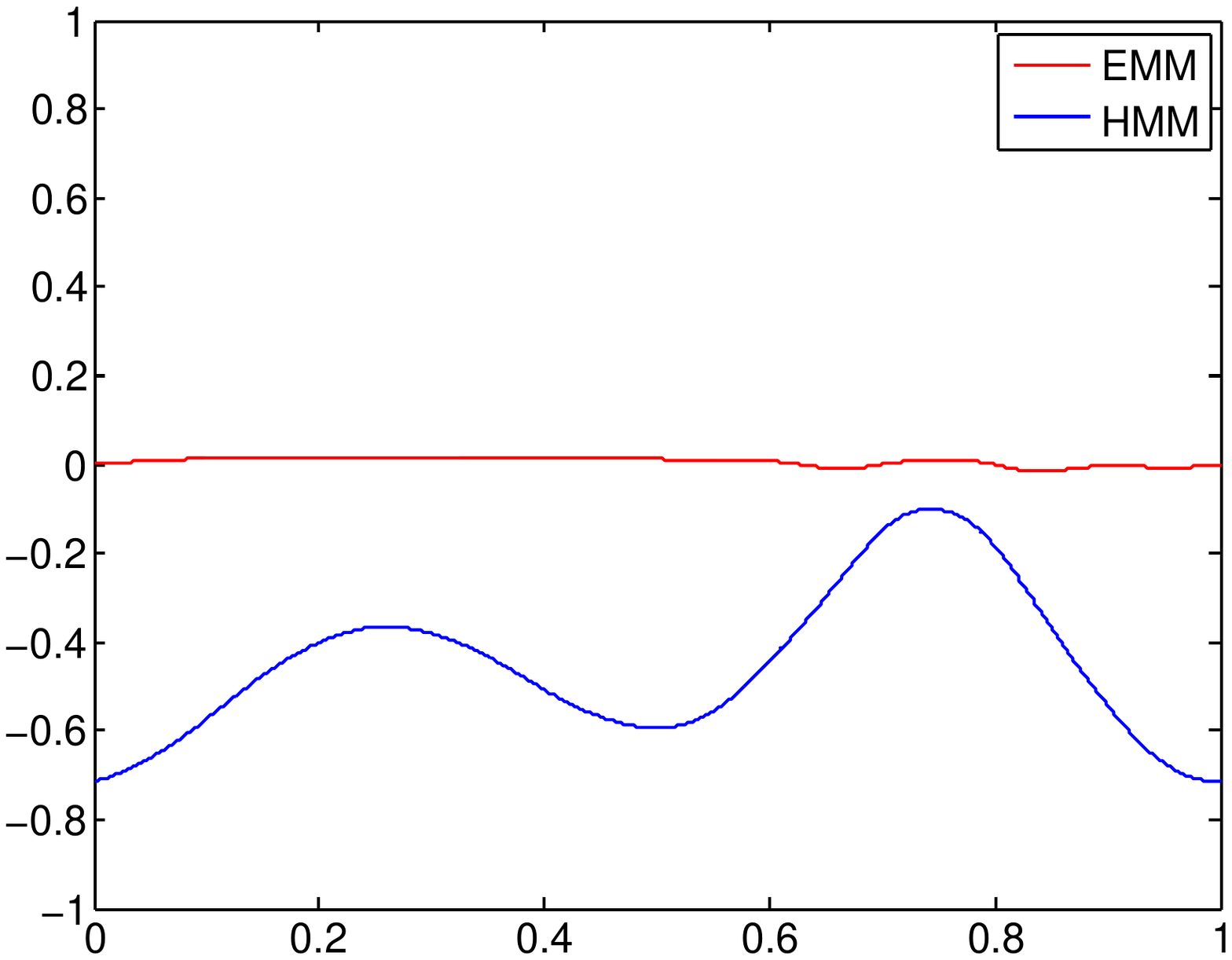} &
\includegraphics[width=5.3cm]{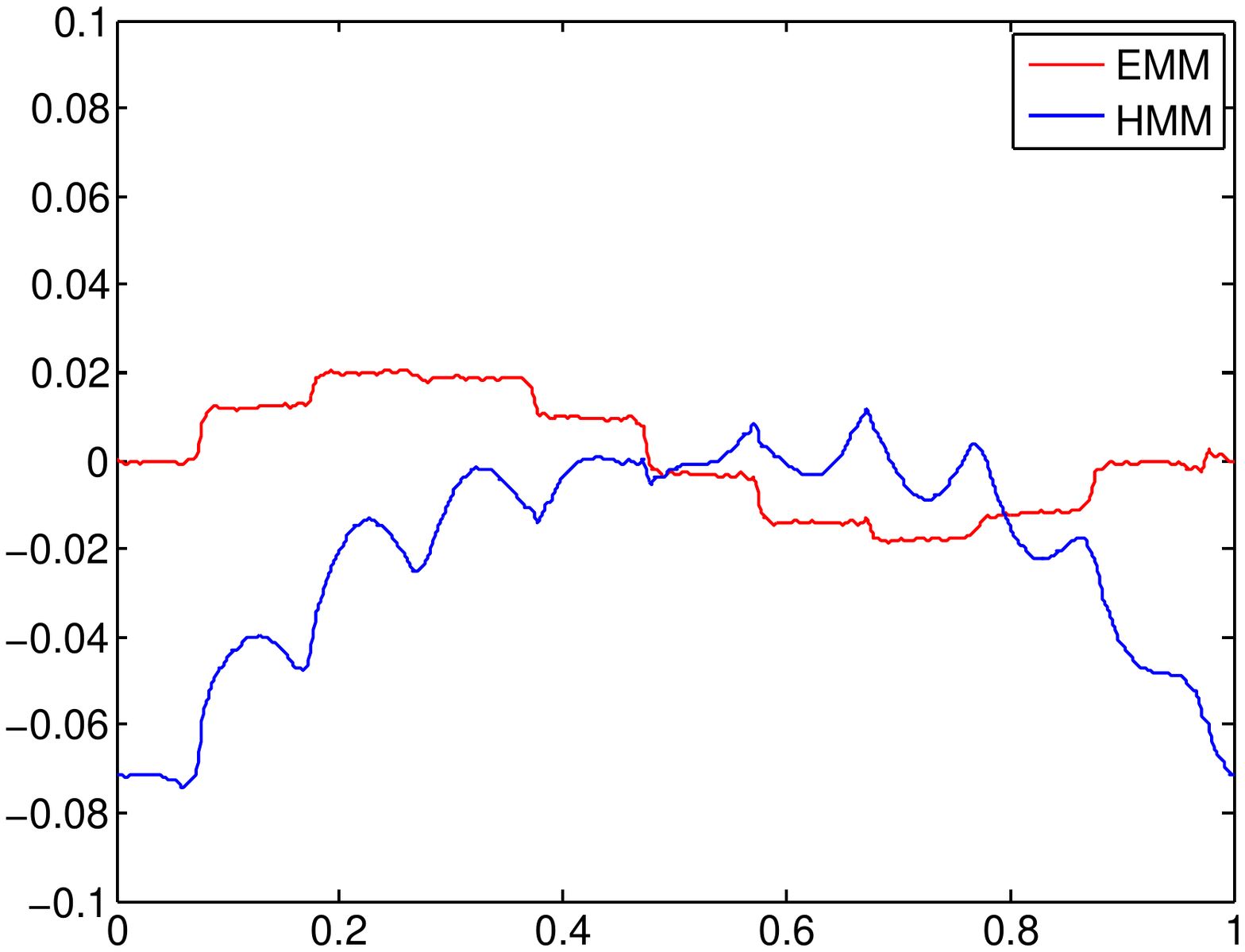} &
\vspace{-4.cm}\includegraphics[width=5.3cm]{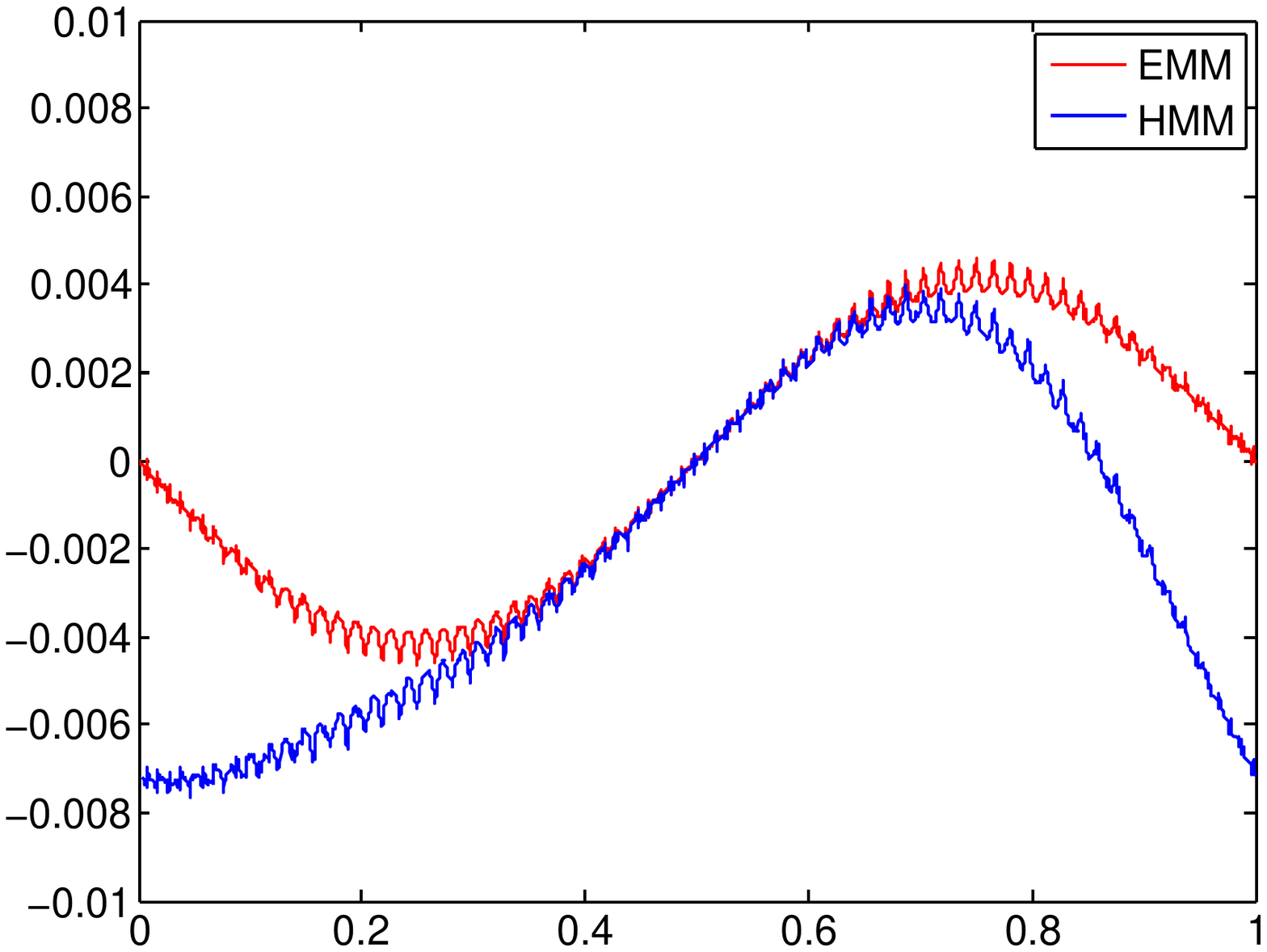} \\
\begin{minipage}{1ex}\vspace{-7.5cm} \rotatebox[origin=lB]{90}{errors for $\partial_x u^\eps$} \end{minipage}&
\includegraphics[width=5.3cm]{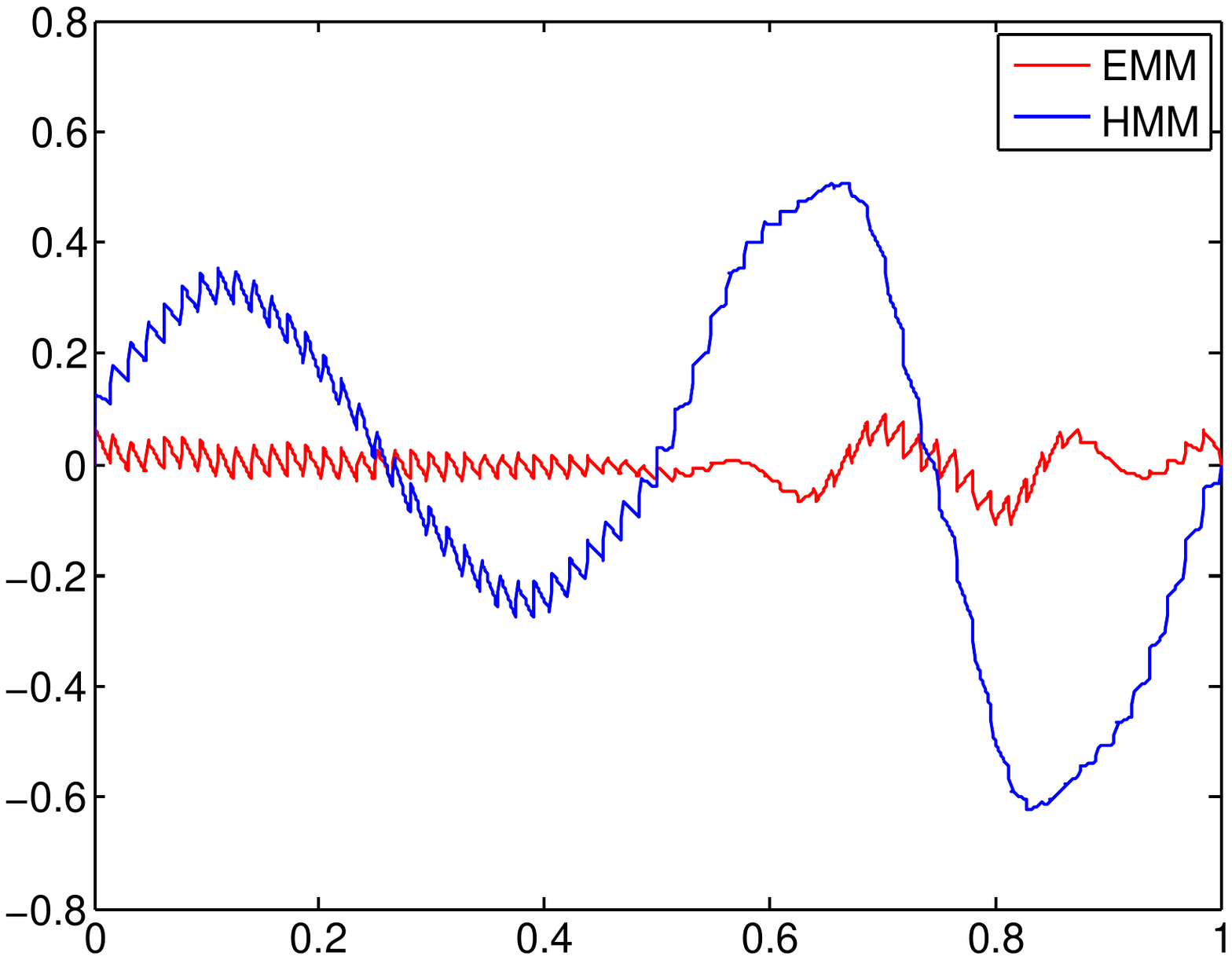} &
\includegraphics[width=5.3cm]{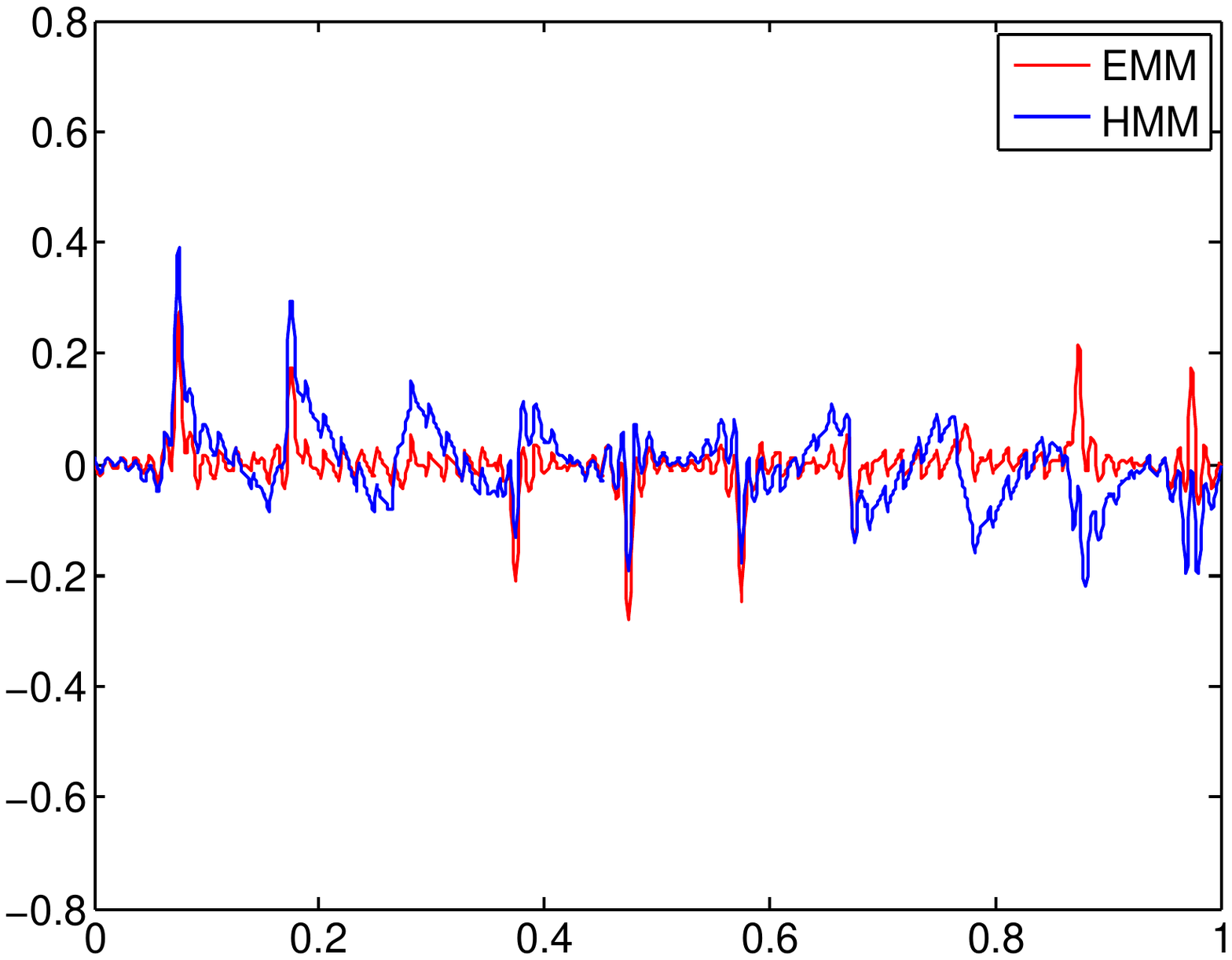} &
\includegraphics[width=5.3cm]{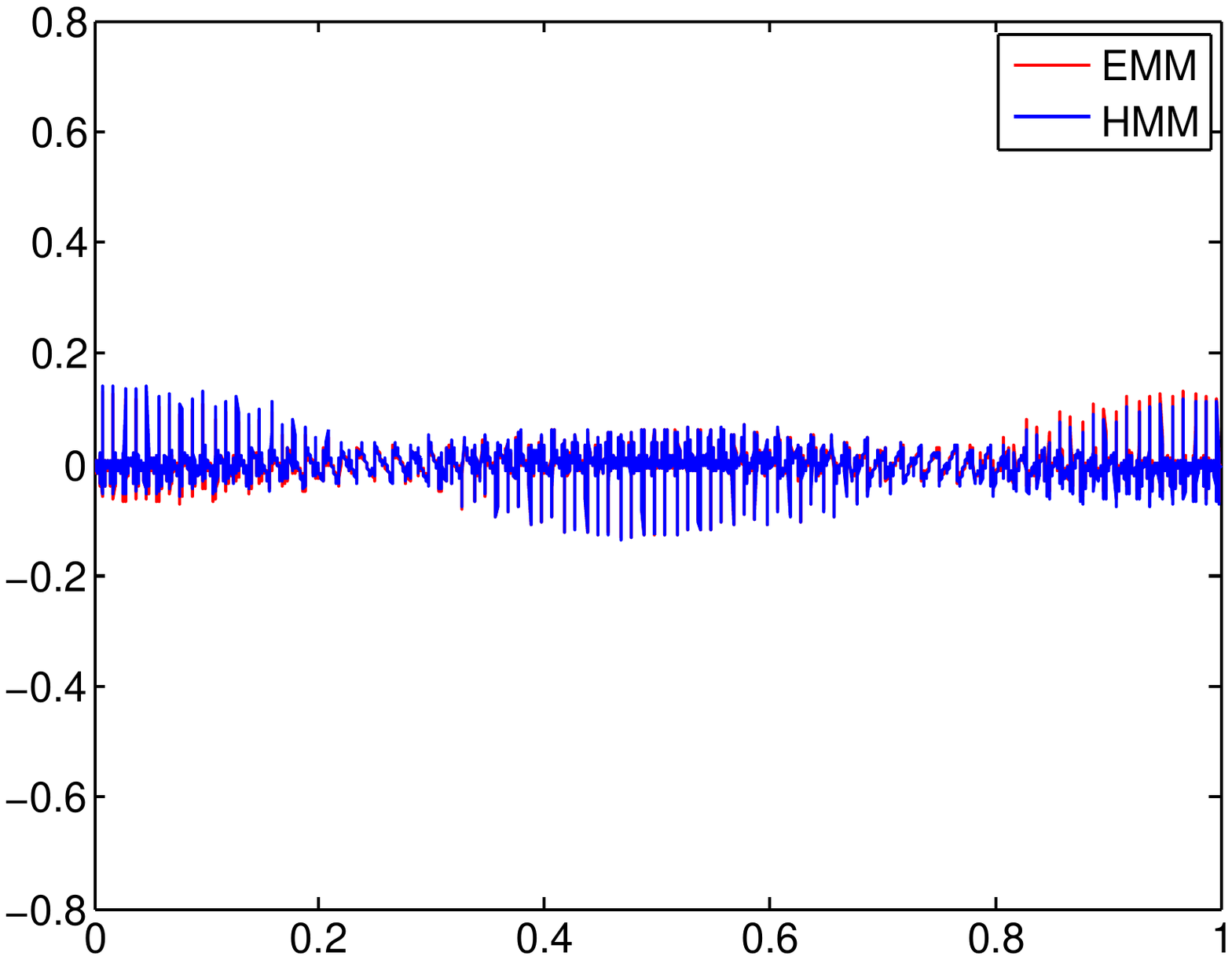} \\
\vspace{-2.8cm}
\end{tabular}
\caption{First line: comparison between $u_{REF}$ (green line), 
$u_{EMM}$ (red line) and $u_{HMM}$ (blue line) for $\varepsilon=1, 0.1, 0.01$ (from left to right). 
Second line (errors in $u^\eps$): $u_{\mbox{\tiny{REF}}} - u_{\mbox{\tiny{EMM}}}$ (red line) and $u_{\mbox{\tiny{REF}} }- u_{\mbox{\tiny{HMM}}}$ (blue line) 
 for $\varepsilon=1, 0.1, 0.01$  (from left to right). 
 Third line (errors in $\partial_x u^\eps$): $\partial_x u_{\mbox{\tiny{REF}}} - \partial_x u_{\mbox{\tiny{EMM}}}$ (red line) and $\partial_x u_{\mbox{\tiny{REF}} }- \partial_x u_{\mbox{\tiny{HMM}}}$ (blue line) 
 for $\varepsilon=1, 0.1, 0.01$ (from left to right). 
The EMM and HMM solutions and their derivatives are computed using the interpolation \eqref{ureconstruct}.}
\label{fig:1}
\end{figure}

\noindent In Figure \ref{fig:1}, we consider the cases $\varepsilon=1, 0.1, 0.01$, respectively (from left to right).
We plot at the final time $T=1$ the solutions given by REF, EMM and HMM as functions of $x\in(0,1)$ (first line), the error in $u^\eps$ (second line, plotting $u_{\mbox{\tiny{REF}}} - u_{\mbox{\tiny{EMM}}}$ and $u_{\mbox{\tiny{REF}} }- u_{\mbox{\tiny{HMM}}}$), and the error in the derivative $\partial_x u^\eps$ (third line, plotting $\partial_x  u_{\mbox{\tiny{REF}}} - \partial_x  u_{\mbox{\tiny{EMM}}}$ and $\partial_x  u_{\mbox{\tiny{REF}} }- \partial_x  u_{\mbox{\tiny{HMM}}}$).
%$v_{\mbox{\tiny{REF}}} - v_{\mbox{\tiny{EMM}}}$ and $v_{\mbox{\tiny{REF}} }- v_{\mbox{\tiny{HMM}}}$ with $v=u$ 
%(second line) and $v=\partial_x u$ (last line), 
%for $\varepsilon=1, 0.1, 0.01$ (from left to right). 
For the reference solution REF, we use $N_x=1024$ for $\varepsilon=1, 0.1$,
and $N_x=4096$ for $\varepsilon=0.01$. 
Concerning EMM and HMM, we use in all computations the mesh parameters $N_x=64$ and $N_y=16$ 
for all $\varepsilon$.  
We can observe that for arbitrary values of $\varepsilon$, EMM is in a very good agreement with 
the refined REF solution and its derivative, for a {\it fixed} set of numerical parameters. 
However, the HMM solution is accurate only in the case $\varepsilon=0.01$. 
In addition, as $\varepsilon$ goes to zero, the computational cost 
for EMM is constant whereas the one of a direct method such as REF increases as the meshes need to be refined. 

%\vspace{-2.cm}

\vspace{0.5ex}

\noindent
\textit{Acknowledgements.}
N.\,C. and M.\,L. acknowledge support by the ANR project Moonrise (ANR-14-CE23-0007-01).
N.\,C. is also partly supported by the ERC Starting Grant Project GEOPARDI.
The research of G.\,V. is partially supported by the Swiss National Foundation, Grant No\,200020\_144313/1. 

\begin{spacing}{0.6}
\bibliographystyle{abbrv}
\bibliography{biblio}
\end{spacing}

%\section{Materiel (to be removed)}
%Tentative scheme to avoid a CFL due to the explicit Euler method:
%\begin{eqnarray}
%\label{noCFL}
%\overline G &=&(I-\Pi) (I-\Delta t D)^{-1} (G_n + \Delta t(I-\Pi) DF^n) \\
 %G^{n+1} &=&\left(I-\frac{\Delta t}{\varepsilon^2}L\right)^{-1}\left[\overline G 
%+ \frac{\Delta t}{\varepsilon}(I-\Pi)B(F^n+\overline G) 
%\right], \\
 %G^{n+1} &=&\left(I-\frac{\Delta t}{\varepsilon^2}L\right)^{-1}(I-\Pi)\left[(I+ \frac{\Delta t}{\varepsilon}B)
%\overline G 
%+ \frac{\Delta t}{\varepsilon}BF^n
%\right], \\
%F^{n+1} &=& \left(I- \Delta t(1-e^{-\Delta t/\varepsilon^2}) \overline{D}\right)^{-1} \left[ F^n+\Delta t (e^{-\Delta t/\varepsilon^2}\Pi DF^n+\frac{  e^{-\Delta t/\varepsilon^2} }{\varepsilon}\Pi B  G^n +  \Pi DG^{n+1}+  S)\right]. 
%\end{eqnarray}

\end{document}